\documentclass{article}
\def\nkyear{201?}
\def\nkvol{3?}
\def\nkno{?}

\def\firstpage{1}
\def\lastpage{20}
\def\correction{4}
\def\shorttitle{Semi-implicit method of sphere-constrained saddle dynamics} 
\def\shortauthor{{\it L. Zhang, P. Zhang\ \ and\ \ X. Zheng}} 
\usepackage b \baselineskip 14.5pt

\input{cam.doi}
\journalnumber{11401}
\DOImsnr{0001}
 \DOIyear{2007}

\title{\Large \bf \boldmath\ \\ Error estimate for semi-implicit method of sphere-constrained high-index saddle dynamics$^{\ast}$} 

\author{\large  Lei ZHANG$^1$ \qquad Pingwen ZHANG$^{2,3}$\qquad Xiangcheng ZHENG$^{4}$} 

\date{}

\begin{document}

\maketitle

\thispagestyle{first}
\renewcommand{\thefootnote}{\fnsymbol{footnote}}

\footnotetext{\hspace*{-5mm} \begin{tabular}{@{}r@{}p{13.4cm}@{}}
& Manuscript received  \\ 
$^1$ & Beijing International Center for Mathematical Research, Center for Machine Learning Research, Center for Quantitative Biology, Peking University, Beijing, 100871, China.\\
&{E-mail:} zhangl@math.pku.edu.cn \\
$^{2}$ & School of Mathematics and Statistics, Wuhan University, Wuhan, 430072, China.\\
$^{3}$ & School of Mathematical Sciences, Laboratory of Mathematics and Applied Mathematics, Peking University, Beijing, 100871, China.\\
&{E-mail:} pzhang@pku.edu.cn \\
$^{4}$ & School of Mathematics, Shandong University, Jinan, 250100, China.\\
&{E-mail:} xzheng@sdu.edu.cn\\
 $^{\ast}$ & This work was partially supported by the National Key R\&D Program of China No. 2021YFF1200500 and the National Natural Science Foundation of China No. 12225102, 12050002, 12288101.
\end{tabular}}

\renewcommand{\thefootnote}{\arabic{footnote}}

\begin{abstract} 
We prove error estimates for the semi-implicit numerical scheme of sphere-constrained high-index saddle dynamics, which serves as a powerful instrument in finding saddle points and constructing the solution landscapes of constrained systems on the high-dimensional sphere. Due to the semi-implicit treatment and the novel computational procedure, the orthonormality of numerical solutions at each time step could not be fully employed to simplify the derivations, and the computations of the state variable and directional vectors are coupled with the retraction, the vector transport and the orthonormalization procedure, which significantly complicates the analysis. We address these issues to prove error estimates for the proposed semi-implicit scheme and then carry out numerical experiments to substantiate the theoretical findings.

\vskip 4.5mm

\nd \begin{tabular}{@{}l@{ }p{10.1cm}} {\bf Keywords } &
Saddle point, Constrained saddle dynamics, Solution landscape, Semi-implicit, Numerical analysis
\end{tabular}

\nd {\bf 2000 MR Subject Classification } 
37M05, 37N30, 65L20

\end{abstract}

\baselineskip 14pt

\setlength{\parindent}{1.5em}

\setcounter{section}{0}

\section{Introduction}
 High-index saddle dynamics \cite{YinSISC} attracts increasing attentions in the last few years due to its capability of effectively finding multiple high-index saddle points of complex systems \cite{EZho,ZhaDu,npj2016}. Here the index of saddle point refers to the Morse index characterized by the maximal dimension of a subspace on which its Hessian operator is negative definite \cite{Milnor}. In particular, the high-index saddle dynamics could be further combined with the downward and upward algorithms \cite{YinSCM} to construct the solution landscape, the pathway map consisting of all stationary points and their connections \cite{wang2021modeling}, that arises several successful applications \cite{Han2019transition,HanXu,YinPRL,Yin2020nucleation,YuZhe, ZhangChe,ZhangCheDu}. In practical problems such as the Thomson problem \cite{Tho} and the Bose-Einstein Condensation \cite{bao2013mathematical}, the state variable is constrained on a high-dimensional sphere, which leads to the more complicated sphere-constrained high-index saddle dynamics for treating the sphere-constrained problems.  

There exist extensive works about numerical analysis to algorithms of finding index-1 saddle points \cite{baker1986,Doye,EV2010,Farr,Gao,Gou,Gra,Lev,Li2001,Mehta,Xie}, while the corresponding analysis for high-index saddle point searchers is rare. In \cite{Z3}, an explicit scheme for the unconstrained high-index saddle dynamics was rigorously analyzed by overcoming the difficulties caused by the coupling of solutions and the (nonlinear) orthonormal procedure of directional vectors in the numerical scheme. The developed method was then extended to prove error estimates for the explicit scheme of the sphere-constrained high-index saddle dynamics by  accounting for the more complex dynamical form and additional operations in the numerical scheme such as the retraction and vector transport in order to maintain the manifold constraint \cite{CHiSD2021}. To improve the numerical stability, a semi-implicit numerical scheme for the unconstrained high-index saddle dynamics was recently analyzed in \cite{LuoSI}, and various numerical experiments demonstrated that comparing with the explicit scheme, the semi-implicit method could improve the convergence behavior, admit much larger step size and reduce the number of queries for the model.

The current work is a continuation of the aforementioned sequence of investigations for numerical analysis of high-index saddle dynamics, which will develop and analyze the semi-implicit numerical method for the sphere-constrained high-index saddle dynamics. To achieve this goal, not only do we need to accommodate the complicated nonlinear forms of this dynamical system, the retraction of the state variable, the vector transport and orthonormalization of the directional vectors due to the manifold constraint, but novel techniques are required to overcome the difficulties caused by the semi-implicit treatment.   The derived results provide theoretical supports for the numerical accuracy of discretization of sphere-constrained high-index saddle dynamics and construction of solution landscapes for complex systems. 

The rest of the paper is organized as follows: In Section 2 we present formulations of the sphere-constrained high-index saddle dynamics and its semi-implicit numerical scheme. In Section 3 we prove several auxiliary estimates, based on which we derive error estimates for the semi-implicit scheme of sphere-constrained high-index saddle dynamics in Section 4. Numerical experiments are performed in Section 5 to substantiate the theoretical findings, and we address concluding remarks in the last section.

\section{Problem formulation and semi-implicit scheme}
In this section we propose the semi-implicit numerical scheme of the sphere-constrained high-index saddle dynamics. Let $E(x)$ be the energy function with $x\in\mathbb R^d$, and define $F(x)=-\nabla E(x)$ and $J(x)=-\nabla^2 E(x)$ with $J(x)=J(x)^\top$. The  high-index saddle dynamics for an index-k saddle point of $E(x)$ constrained on the unit sphere $S^{d-1}$ was developed in \cite{CHiSD2021}:
\begin{equation}\label{csadk}
\left\{
\begin{array}{l}
\ds \frac{dx}{dt} =\bigg(I -xx^\top-2\sum_{j=1}^k v_jv_j^\top \bigg)F(x);\\[0.075in]
\ds \frac{dv_i}{dt}=\bigg( I-xx^\top-v_iv_i^\top-2\sum_{j=1}^{i-1}v_jv_j^\top\bigg)J(x)v_i+ xv_i^\top F(x)
\end{array}
\right.
\end{equation}
for $1\leq i\leq k$, equipped with the initial conditions
\begin{align*}
&x(0)=x_0\in S^{d-1},~~v_i(0)=v_{i,0}\\
\text{ such that }&v_{i,0}^\top v_{j,0}=\delta_{ij}\text{ and }x_0^\top v_{i,0}=0\text{ for }1\leq i,j\leq k. 
\end{align*}
Here $x$ represents a position variable and $\{v_i\}_{i=1}^k$ are $k$ directional variables. It was proved in \cite{CHiSD2021} that a linearly stable steady state of (\ref{csadk}) is an index-$k$ saddle point, and  the solutions $x$ and $\{v_i\}_{i=1}^k$ to the dynamics (\ref{csadk}) satisfy for $t>0$
\begin{equation}\label{prop}
 x(t)\in S^{d-1},~~v_i(t)^\top x(t)=0,~~v_{i}(t)^\top v_{j}(t)=\delta_{ij},~~1\leq i,j\leq k.
 \end{equation}

Throughout the paper we apply the following assumptions: 

\noindent\textbf{Assumption A:} The $F(x)$ could be represented as a sum of the linear part $\mathcal Lx $ and the nonlinear part $\mathcal N(x)$, that is, $F(x)=\mathcal Lx+\mathcal N(x)$, and there exists a constant $L>0$ such that the following linearly growth and Lipschitz conditions hold under the standard $l^2$ norm $\|\cdot\|$ of a vector or a matrix
$$\begin{array}{c}
\ds \max\{\|J(x_2)-J(x_1)\|,\|\mathcal Lx_2-\mathcal Lx_1\|,\|\mathcal N(x_2)-\mathcal N(x_1)\|\}\leq L\|x_2-x_1\|,\\[0.1in]
\ds\max\{\|\mathcal Lx\|,\|\mathcal N(x)\|\}\leq L(1+\|x\|),~~x,x_1,x_2\in \mathbb R^d.
\end{array}  $$

To derive the semi-implicit discretization, let $0=t_0<t_1<\cdots t_N=T$ be the uniform partition of $[0,T]$ with the step size $\tau=T/N$, and let $\{x_n,v_{i,n}\}_{n=0}^N$ be the numerical solution of (\ref{csadk}). Then we discretize the first-order derivative by the Euler scheme and treat the linear and nonlinear parts on the right-hand side of (\ref{csadk}) via the implicit and explicit manner, respectively, to obtain the semi-implicit scheme of (\ref{csadk}) for $1\leq n\leq N$ as follows:
\begin{equation}\label{cFDsadk}
\left\{
\begin{array}{l}
\ds \tilde x_{n} =x_{n-1}+\tau\bigg(I -2\sum_{j=1}^k v_{j,n-1}v_{j,n-1}^\top \bigg)(\mathcal L \tilde x_{n}+\mathcal N(x_{n-1}))\\[0.2in]
 \qquad\qquad\ds-\tau x_{n-1}x_{n-1}^\top(\mathcal L x_{n-1}+\mathcal N(x_{n-1})),\\[0.1in]
\ds x_n=\frac{\tilde x_n}{\|\tilde x_n\|};\\[0.2in]
\left.
\begin{array}{l}
\ds \tilde v_{i,n}=v_{i,n-1}+\tau\bigg( I-x_{n}x_{n}^\top-2\sum_{j=1}^{i-1}v_{j,n}v_{j,n}^\top\bigg)J(x_{n})\tilde v_{i,n}\\[0.2in]
\ds \qquad\qquad-\tau v_{i,n-1}v_{i,n-1}^\top J(x_{n})v_{i,n-1}+\tau x_{n}\tilde v_{i,n}^\top F(x_{n}),\\[0.1in]
\hat v_{i,n}=\tilde v_{i,n}-\tilde v_{i,n}^\top x_n x_n,\\[0.05in]
v_{i,n}=\mbox{GS}(\hat v_{i,n},\{v_{j,n}\}_{j=1}^{i-1}),\\
\end{array}
\right\}
~1\leq i\leq k.
\end{array}
\right.
\end{equation}
Here the Gram-Schmidt orthonormalization function GS($\hat v_{i,n},\{v_{j,n}\}_{j=1}^{i-1}$) generates the normalized vector $v_{i,n}$ from $\hat v_{i,n}$ that is orthogonal with $\{v_{j,n}\}_{j=1}^{i-1}$, that is,
$$ 
\ds  v_{i,n}=\mathcal N \bigg(\ds\hat v_{i,n}-\sum_{j=1}^{i-1}(\hat v_{i,n}^\top v_{j,n})v_{j,n}\bigg):=\frac{1}{Y_{i,n}}\bigg(\ds\hat v_{i,n}-\sum_{j=1}^{i-1}(\hat v_{i,n}^\top v_{j,n})v_{j,n}\bigg),$$
where $\mathcal N$ is the normalized operator  and the normalized factor $Y_{i,n}$ is thus defined as
$$
\ds Y_{i,n}:=\bigg\|\hat v_{i,n}-\sum_{j=1}^{i-1}(\hat v_{i,n}^\top v_{j,n})v_{j,n}\bigg\|=\bigg(\|\hat v_{i,n}\|^2-\sum_{j=1}^{i-1}(\hat v_{i,n}^\top v_{j,n})^2\bigg)^{1/2}.
$$
The first and the third schemes in (\ref{cFDsadk}) are semi-implicit discretizations of the equations of $x$ and $v_i$ in (\ref{csadk}), respectively. The second equation of (\ref{cFDsadk}) represents the retraction in order to ensure that $x_n\in S^{d-1}$. The last two schemes, which stand for the vector transport and the Gram-Schmidt orthonormalization procedure, respectively, aim to ensure the rest properties of (\ref{prop}), that is,
\begin{equation}\label{propn}
v_{i,n}^\top x_n=0,~~v_{i,n}^\top v_{j,n}=\delta_{ij},~~1\leq i,j\leq k,~~0\leq n\leq N.
\end{equation}

Different from the explicit scheme presented in \cite{Z3c}, where all variables on the right-hand side of (\ref{cFDsadk}) take their values at the previous time step $t_{n-1}$, the orthonormal property of the vectors $\{v_{i,n-1}\}_{i=1}^k$ at the time step $t_{n-1}$ could no longer be fully employed in (\ref{cFDsadk}) to facilitate the numerical analysis as performed in \cite{Z3c} due to the semi-implicit treatment, which complicates the error estimate. On the other hand, in the explicit scheme the vectors $\{\tilde v_{i,n}\}_{i=1}^k$ are firstly solved, and then their orthonormalization are independently performed. In the semi-implicit scheme (\ref{cFDsadk}), the computational strategy is quite different in that the last three schemes of directional vectors in (\ref{cFDsadk}) are sequentially solved for $1\leq i\leq k$. In this way, the newly computed orthonormalized vectors $\{v_{j,n}\}_{j=1}^{i-1}$ at the current time step $t_n$ are involved in the scheme of $\tilde v_{i,n}$, which could be more appropriate than invoking the vectors at the previous time step in the explicit scheme. However, this computational strategy leads to the coupling of the schemes of directional vectors, the vector transport and the orthonormalization procedure, which makes the numerical analysis more challenging.

Concerning these difficulties, we derive novel analysis methods to carry out error estimates in subsequent sections. Throughout the paper we use $Q$ to denote a generic positive constant that may assume different values at different occurrences.

\section{Auxiliary estimates}
We prove several properties of the numerical solutions to support the error estimates. By $\|x_n\|=\|v_{i,n}\|=1$ for $1\leq i\leq k$ and $1\leq n\leq N$, we could apply the Assumption A to derive from the first and the third equations of the scheme (\ref{cFDsadk}) that 
\begin{equation}\label{xbndmh}
\max\{\|\tilde x_n\|,\|\tilde v_{1,n}\|,\cdots,\|\tilde v_{k,n}\|\}\leq Q
\end{equation}
 for $1\leq n\leq N$ for $\tau$ small enough, which will be frequently used in the analysis.
\begin{lem}\label{clem1}
Under the Assumption A, the following estimate holds for $\tau$ small enough:
\begin{align}
\label{xtx}
&\|x_n-\tilde x_n\|\leq Q\tau^2,~~1\leq n\leq N;\\
&\|\hat v_{i,n}-\tilde v_{i,n}\|=|\tilde v_{i,n}^\top x_n|\leq Q\tau^2,~~1\leq i\leq k,~~1\leq n\leq N. \label{vmh}
\end{align}
\end{lem}
\begin{proof}
We employ the first equation of (\ref{cFDsadk}) to get
\begin{equation}\label{xxmh}
\begin{array}{l}
\ds\|\tilde x_n-x_{n-1}\|=\bigg\|\tau\bigg(I -2\sum_{j=1}^k v_{j,n-1}v_{j,n-1}^\top \bigg)(\mathcal L \tilde x_{n}+\mathcal N(x_{n-1}))\\[0.2in]
 \qquad\qquad\qquad\qquad\ds-\tau x_{n-1}x_{n-1}^\top(\mathcal L x_{n-1}+\mathcal N(x_{n-1}))\bigg\|\leq Q\tau.
\end{array} 
\end{equation}
We then apply this to rewrite the first equation of (\ref{cFDsadk}) as
\begin{align}
\tilde x_{n} &=x_{n-1}+\tau\bigg(I -2\sum_{j=1}^k v_{j,n-1}v_{j,n-1}^\top \bigg)(\mathcal L \tilde x_{n}+\mathcal N(x_{n-1}))\nonumber\\
 &\qquad-\tau x_{n-1}x_{n-1}^\top(\mathcal L x_{n-1}+\mathcal N(x_{n-1}))\nonumber\\
 &=x_{n-1}+\tau\bigg(I- x_{n-1}x_{n-1}^\top-2\sum_{j=1}^k v_{j,n-1}v_{j,n-1}^\top \bigg)(\mathcal L \tilde x_{n}+\mathcal N(x_{n-1}))\label{lalala}\\
 &\qquad+\tau x_{n-1}x_{n-1}^\top\mathcal L(\tilde x_n-x_{n-1})\nonumber\\
 &=x_{n-1}+\tau\bigg(I- x_{n-1}x_{n-1}^\top-2\sum_{j=1}^k v_{j,n-1}v_{j,n-1}^\top \bigg)(\mathcal L \tilde x_{n}+\mathcal N(x_{n-1}))+O(\tau^2).\nonumber
\end{align}
We multiply $x_{n-1}^\top$ on both sides of this equation and use (\ref{propn}) to obtain 
$$x_{n-1}^\top\tilde x_n=1+O(\tau^2).$$
 We then multiply $\tilde x_{n}^\top$ on both sides of (\ref{lalala}) and use $x_{n-1}^\top v_{j,n-1}=0$ for $1\leq j\leq k$ and $x_{n-1}^\top\tilde x_n=1+O(\tau^2)$ to obtain
\begin{align*}
\|\tilde x_{n}\|^2&=1+\tau\bigg(\tilde x_n^\top-x_{n-1}^\top -2\sum_{j=1}^k\tilde x_n^\top v_{j,n-1}v_{j,n-1}^\top \bigg)(\mathcal L \tilde x_{n}+\mathcal N(x_{n-1}))+O(\tau^2)\\
&= 1+\tau(\tilde x_n-x_{n-1})^\top\bigg(I-2\sum_{j=1}^kv_{j,n-1}v_{j,n-1}^\top\bigg)(\mathcal L \tilde x_{n}+\mathcal N(x_{n-1}))+O(\tau^2),
\end{align*}
which, together with the Assumption A and the norm-preserving property of the Householder matrix in the above equation,
 yields
\begin{align*}
\big|\|\tilde x_n\|^2-1\big|&\leq \tau\|\tilde x_n-x_{n-1}\|\|\mathcal L \tilde x_{n}+\mathcal N(x_{n-1})\|+O(\tau^2)\\
&\leq Q\tau\|\tilde x_n-x_{n-1}\|+O(\tau^2). 
\end{align*}
Combining this equation and (\ref{xxmh}) we obtain
$$ \big|\|\tilde x_n\|^2-1\big|\leq Q \tau^2, $$
 which in turn leads to  $ \big|\|\tilde x_n\|-1\big|\leq Q \tau^2 $. We apply this to reach  (\ref{xtx}):
$$\|x_n-\tilde x_n\|=\bigg\|\frac{\tilde x_n}{\|\tilde x_n\|}(1-\|\tilde x_n\|)\bigg\|= \big|1-\|\tilde x_n\|\big|\leq  Q\tau^2.$$

To derive (\ref{vmh}), we combine (\ref{xtx}) and (\ref{xxmh}) to obtain
\begin{equation}\label{xnn}
\|x_n-x_{n-1}\|\leq \|x_n-\tilde x_n\|+\|\tilde x_n-x_{n-1}\|\leq Q\tau.
\end{equation}
From the forth equation of (\ref{cFDsadk}) we apply $\|x_n\|=1$ to obtain
\begin{equation}\label{clem11}
\|\hat v_{i,n}-\tilde v_{i,n}\|=|\tilde v_{i,n}^\top x_n|.
 \end{equation}
Furthermore, the relation $ \big|\|\tilde x_n\|-1\big|\leq Q \tau^2 $ leads to $\|\tilde x_n\|\geq 1-Q\tau^2\geq 1/2$ for $\tau$ small enough.
Then we multiply the scheme of $\tilde v_{i,n}$ in (\ref{cFDsadk}) and the reformulated scheme of $\tilde x_n$ in (\ref{lalala}) to get
\begin{align}
x_n^\top \tilde v_{i,n} &=\frac{1}{\|\tilde x_n\|}\bigg[x_{n-1}+\tau\bigg(I- x_{n-1}x_{n-1}^\top-2\sum_{j=1}^k v_{j,n-1}v_{j,n-1}^\top \bigg)\nonumber\\
&\qquad\cdot(\mathcal L \tilde x_{n}+\mathcal N(x_{n-1}))+O(\tau^2)\bigg]^\top\nonumber\\
 &~~\bigg[ v_{i,n-1}+\tau\bigg( I-x_{n}x_{n}^\top-2\sum_{j=1}^{i-1}v_{j,n}v_{j,n}^\top\bigg)J(x_{n})\tilde v_{i,n}\nonumber\\
&\qquad-\tau v_{i,n-1}v_{i,n-1}^\top J(x_{n})v_{i,n-1}+\tau x_{n}\tilde v_{i,n}^\top F(x_{n})\bigg]\nonumber\\
 &=\frac{\tau}{\|\tilde x_n\|}\bigg[\bigg(x_{n-1}^\top -x_{n-1}^\top x_n x_n^\top-2\sum_{j=1}^{i-1}x_{n-1}^\top v_{j,n}v_{j,n}^\top\bigg)J(x_n)\tilde v_{i,n}\nonumber\\
&\qquad+x_{n-1}^\top x_n \tilde v_{i,n}^\top F(x_n) -v_{i,n-1}^\top (\mathcal L \tilde x_{n}+\mathcal N(x_{n-1}))\bigg]+O(\tau^2),\label{vxmh}
\end{align}
where we briefly write the second-order terms of $\tau$ as $O(\tau^2)$. We apply the splittings
$$x_{n-1}^\top -x_{n-1}^\top x_n x_n^\top=(x_{n-1}-x_n)^\top(I -x_n x_n^\top) $$
and
$$\begin{array}{l}
\ds x_{n-1}^\top x_n \tilde v_{i,n}^\top F(x_n) -v_{i,n-1}^\top (\mathcal L \tilde x_{n}+\mathcal N(x_{n-1}))\\[0.05in]
\ds \qquad=x_{n-1}^\top x_n \tilde v_{i,n}^\top F(x_n) -v_{i,n-1}^\top F(x_{n-1})+v_{i,n-1}^\top \mathcal L (x_{n-1}-\tilde x_{n})\\[0.05in]
\ds\qquad=(x_{n-1}-x_n)^\top x_n \tilde v_{i,n}^\top F(x_n)+(\tilde v_{i,n}-v_{i,n-1})^\top F(x_n)\\[0.05in]
\ds\qquad\qquad+v_{i,n-1}^\top(F(x_n)-F(x_{n-1}))+v_{i,n-1}^\top \mathcal L (x_{n-1}-\tilde x_{n}),
\end{array}  $$
to bound the right-hand side of (\ref{vxmh}) as
\begin{align}
|x_n^\top \tilde v_{i,n}|&\leq \frac{Q\tau}{\|\tilde x_n\|}\Big[\|x_n-x_{n-1}\|+\sum_{j=1}^{i-1}|x_{n-1}^\top v_{j,n}|\nonumber\\
&\qquad+\|\tilde v_{i,n}-v_{i,n-1}\|+\|F(x_n)-F(x_{n-1})\|+\|x_{n-1}-\tilde x_n||\Big]+O(\tau^2).\label{klok}
\end{align}
We then invoke the third scheme of (\ref{cFDsadk})
\begin{align}\nonumber
&\|\tilde v_{i,n}-v_{i,n-1}\|=\bigg\|\tau\bigg( I-x_{n}x_{n}^\top-2\sum_{j=1}^{i-1}v_{j,n}v_{j,n}^\top\bigg)J(x_{n})\tilde v_{i,n}\\
& \qquad\qquad-\tau v_{i,n-1}v_{i,n-1}^\top J(x_{n})v_{i,n-1}+\tau x_{n}\tilde v_{i,n}^\top F(x_{n}) \bigg\|\leq Q\tau,\label{vvmh}
\end{align}
as well as $x_{n-1}^\top v_{j,n}=(x_{n-1}-x_n)^\top v_{j,n}$, $\|\tilde x_n\|\geq 1/2$, (\ref{xxmh}), (\ref{xnn}) and the Lipschitz condition of $F$ in (\ref{vxmh}) to obtain
	$$\begin{array}{l}
\ds |\tilde v_{i,n}^\top x_n|\leq Q\tau\Big[\|x_n-x_{n-1}\|+\sum_{j=1}^{i-1}\|x_{n-1}-x_n\|\|v_{j,n}\|\nonumber\\
\ds\qquad\quad\qquad\qquad+\|\tilde v_{i,n}-v_{i,n-1}\|+\|x_{n-1}-\tilde x_n||\Big]+O(\tau^2)\\[0.05in]
\ds\qquad\qquad\leq Q\tau\big(\|x_{n-1}-x_n\|+\|\tilde v_{i,n}-v_{i,n-1}\|+\|x_{n-1}-\tilde x_n\|\big)+O(\tau^2)\leq Q\tau^2,
\end{array} $$
which completes the proof.
\end{proof}

\begin{lem}\label{clem2}
For $1\leq m<i\leq k$ and $1\leq j\leq k$, the following estimates hold for $\tau$ small enough:
 $$\begin{array}{c}
 \ds \|\tilde v_{i,n}^\top \tilde v_{m,n}\|\leq Q\tau\sum_{l=1}^m\|\hat v_{l,n}-v_{l,n}\|+Q\tau^2,\\
 \ds \big|\|\tilde v_{j,n}\|^2-1\big|\leq Q\tau\sum_{l=1}^{j-1}\|\hat v_{l,n}-v_{l,n}\|+Q\tau^2.
\end{array}   $$
\end{lem}
\begin{proof}
From the definitions of $\tilde v_{i,n}$ and $\tilde v_{m,n}$ we have
$$\begin{array}{l}
\ds \tilde v_{i,n}^\top \tilde v_{m,n}=\tau\bigg(v_{m,n-1}^\top J(x_n)\tilde v_{i,n}-x_n^\top v_{m,n-1}\tilde v_{i,n}J(x_n)^\top x_n\\
\ds\qquad\qquad\quad-2\sum_{j=1}^{i-1}v_{m,n-1}^\top v_{j,n}v_{j,n}^\top J(x_n)\tilde v_{i,n}+ x_n^\top v_{m,n-1}\tilde v_{i,n}^\top F(x_n)\\[0.2in]
\ds\qquad\qquad\quad+\tilde v_{m,n}^\top J(x_n)^\top v_{i,n-1}-x_n^\top v_{i,n-1} x_n^\top J(x_n)\tilde v_{m,n}\\[0.05in]
\ds\qquad\qquad\quad-2\sum_{j=1}^{m-1}v_{j,n}^\top v_{i,n-1} v_{j,n}^\top J(x_n)\tilde v_{m,n}+v_{i,n-1}^\top x_n \tilde v_{m,n}^\top F(x_n)\bigg)+O(\tau^2)\\[0.05in]
\ds\qquad\qquad=:\sum_{l=1}^8 K_l+O(\tau^2).
\end{array} $$
We apply $x_{n}^\top v_{i,n}=0$ for $1\leq i\leq k$ and $1\leq n\leq N$ and (\ref{xnn}) to bound $K_2+K_4+K_6+K_8$ as
$$\begin{array}{l}
\|K_2+K_4+K_6+K_8\|\\[0.05in]
\quad=\tau\big\|-x_n^\top v_{m,n-1}\tilde v_{i,n}J(x_n)^\top x_n+ x_n^\top v_{m,n-1}\tilde v_{i,n}^\top F(x_n)\\[0.05in]
\ds\qquad\quad-x_n^\top v_{i,n-1} x_n^\top J(x_n)\tilde v_{m,n}+v_{i,n-1}^\top x_n \tilde v_{m,n}^\top F(x_n)\big\|\\[0.05in]
\quad=\tau\big\|-(x_n-x_{n-1})^\top v_{m,n-1}\tilde v_{i,n}J(x_n)^\top x_n\\[0.05in]
\ds\quad\qquad+ (x_n-x_{n-1})^\top v_{m,n-1}\tilde v_{i,n}^\top F(x_n)\\[0.05in]
\ds\qquad\quad-(x_n-x_{n-1})^\top v_{i,n-1} x_n^\top J(x_n)\tilde v_{m,n}\\[0.05in]
\ds\quad\qquad+v_{i,n-1}^\top (x_n-x_{n-1}) \tilde v_{m,n}^\top F(x_n)\big\|\leq Q\tau^2.
\end{array}$$
We then introduce the following triple splitting:
\begin{equation*}
v_{i,n-1}-v_{i,n}=(v_{i,n-1}-\tilde v_{i,n})+(\tilde v_{i,n}-\hat v_{i,n})+(\hat v_{i,n}-v_{i,n}).
\end{equation*}
The first right-hand side term is estimated by (\ref{vvmh}) and the second right-hand side term is bounded by Lemma \ref{clem1}, which lead to
\begin{equation}\label{vsplit}
\|v_{i,n-1}-v_{i,n}\|\leq Q\tau+\|\hat v_{i,n}-v_{i,n}\|.
\end{equation}
We invoke this to bound $K_7$ as
$$\begin{array}{rl}
\ds|K_7|&\ds=\bigg|2\tau\gamma\sum_{j=1}^{m-1}v_{j,n}^\top v_{i,n-1} v_{j,n}^\top J(x_n)\tilde v_{m,n}\bigg|\\[0.1in]
&\ds=\bigg|2\tau\gamma\sum_{j=1}^{m-1}(v_{j,n}^\top-v_{j,n-1}^\top) v_{i,n-1} v_{j,n}^\top J(x_n)\tilde v_{m,n}\bigg|\\
&\ds\leq Q\tau^2+Q\tau\sum_{j=1}^{m-1}\|v_{j,n}-\hat v_{j,n}\|.
\end{array}  $$
By $v_{m,n}^\top v_{j,n}=\delta_{m,j}$ we rewrite $K_3$ as
\begin{equation}\label{conc1} \begin{array}{l}
\ds K_3=-2\tau\sum_{j=1}^{i-1}v_{m,n-1}^\top v_{j,n}v_{j,n}^\top J(x_n)\tilde v_{i,n}\\
\ds\quad~
=-2\tau\sum_{j=1}^{i-1}(v_{m,n-1}^\top-v_{m,n}^\top) v_{j,n}v_{j,n}^\top J(x_n)\tilde v_{i,n}-2\tau v_{m,n}^\top J(x_n) \tilde v_{i,n},
\end{array}  \end{equation}
which leads to
\begin{equation}\label{conc2} \begin{array}{rl}
\ds K_1+K_3+K_5&\ds=\tau\big(v_{m,n-1}^\top J(x_n)\tilde v_{i,n}-v_{m,n}^\top J(x_n) \tilde v_{i,n} \big)\\[0.1in]
&\ds\hspace{-0.5in}+\tau\big(\tilde v_{m,n}^\top J(x_n)^\top v_{i,n-1}-v_{m,n}^\top J(x_n) \tilde v_{i,n}\big)\\[0.1in]
&\ds\hspace{-0.5in}-2\tau\sum_{j=1}^{i-1}(v_{m,n-1}^\top-v_{m,n}^\top) v_{j,n}v_{j,n}^\top J(x_n)\tilde v_{i,n}=:B_1+B_2+B_3.
\end{array} \end{equation}

We then use (\ref{vsplit}) to bound $B_1$ as
$$|B_1|=\tau|(v_{m,n-1}^\top -v_{m,n}^\top) J(x_n) \tilde v_{i,n}|\leq Q\tau^2+Q\tau \|\hat v_{m,n}-v_{m,n}\|.$$
 $B_3$ could be estimated similarly:
$$\begin{array}{l}
\ds |B_3|=2\tau\bigg|\sum_{j=1}^{i-1}(v_{m,n-1}^\top-v_{m,n}^\top) v_{j,n}v_{j,n}^\top J(x_n)\tilde v_{i,n}\bigg|\leq Q\tau^2+Q\tau \|\hat v_{m,n}-v_{m,n}\|.
\end{array} $$
We then apply the symmetry of $J(x_n)$ and Lemma \ref{clem1} and (\ref{vvmh})
 to bound $B_2$ as
\begin{align*}
\ds|B_2|&=\tau|(\tilde v_{m,n}^\top-v_{m,n}^\top) J(x_n) v_{i,n-1}+v_{m,n}^\top J(x_n) (v_{i,n-1}-\tilde v_{i,n})|\\
&=\tau|(\tilde v_{m,n}^\top-\hat v_{m,n}^\top+\hat v_{m,n}^\top-v_{m,n}^\top) J(x_n) v_{i,n-1}\\
&\quad+v_{m,n}^\top J(x_n) (v_{i,n-1}-\tilde v_{i,n})| \leq Q\tau^2+Q\tau \|\hat v_{m,n}-v_{m,n}\|.
\end{align*}
We incorporate the preceding estimates to complete the proof of the first statement of this lemma.

To derive the second statement, we apply the definition of $\tilde v_{j,n}$ in (\ref{cFDsadk}) to get
$$\begin{array}{l}
\ds \|\tilde v_{j,n}\|^2=1+2\tau\bigg(v_{j,n-1}^\top-v_{j,n-1}^\top x_n x_n^\top-2\sum_{l=1}^{j-1}v_{j,n-1}^\top v_{l,n} v_{l,n}^\top\bigg)J(x_n)\tilde v_{j,n}\\[0.2in]
\ds\qquad\qquad\qquad-2\tau v_{j,n-1}^\top J(x_n)v_{j,n-1}+2\tau v_{j,n-1}^\top x_n \tilde v_{j,n}^\top F(x_n)+O(\tau^2),
\end{array}  $$
that is,
$$\begin{array}{l}
\ds \big|\|\tilde v_{j,n}\|^2-1\big|=\bigg|2\tau v_{j,n-1}^\top J(x_n)(\tilde v_{j,n}-v_{j,n-1})\\[0.1in]
\ds\qquad\qquad\qquad\quad- \tau\bigg(v_{j,n-1}^\top (x_n-x_{n-1}) x_n^\top\\[0.1in]
\ds\qquad\qquad\qquad\quad+2\sum_{l=1}^{j-1}v_{j,n-1}^\top (v_{l,n}-v_{l,n-1}) v_{l,n}^\top\bigg)J(x_n)\tilde v_{j,n}\\[0.1in]
\ds\qquad\qquad\qquad\quad+2\tau v_{j,n-1}^\top (x_n-x_{n-1}) \tilde v_{j,n}^\top F(x_n)+O(\tau^2)\bigg|.
\end{array}  $$
Thus we incorporate (\ref{xnn}), (\ref{vvmh})
and (\ref{vsplit}) to get 
$$\begin{array}{l}
\ds \big|\|\tilde v_{j,n}\|^2-1\big|\leq Q\tau\sum_{l=1}^{j-1}\|\hat v_{l,n}-v_{l,n}\|+Q\tau^2,
\end{array}  $$
 which completes the proof.
\end{proof}
\begin{lem}\label{clem3}
For $1\leq m<i\leq k$ and $1\leq j\leq k$, the following estimates hold for $\tau$ small enough:
 $$\begin{array}{c}
 \ds \|\hat v_{i,n}^\top \hat v_{m,n}\|\leq Q_0\tau\sum_{l=1}^m\|\hat v_{l,n}-v_{l,n}\|+Q_1\tau^2,\\
 \ds \big|\|\hat v_{j,n}\|^2-1\big|\leq Q_2\tau\sum_{l=1}^{j-1}\|\hat v_{l,n}-v_{l,n}\|+Q_3\tau^2.
\end{array}   $$
\end{lem}
\begin{proof}
 For $1\leq m< i\leq k$ we get
 $$\hat v_{m,n}^\top \hat v_{i,n}=\tilde v_{m,n}^\top \tilde v_{i,n}-x_n^\top \tilde v_{i,n}x_n^\top \tilde v_{m,n}, $$
which, together with Lemmas \ref{clem1} and \ref{clem2}, leads to
\begin{align*}
\big|\hat v_{m,n}^\top \hat v_{i,n}\big|&\leq Q\tau\sum_{l=1}^m\|\hat v_{l,n}-v_{l,n}\|+Q\tau^2+Q\tau^4\\
&\leq Q\tau\sum_{l=1}^m\|\hat v_{l,n}-v_{l,n}\|+Q\tau^2.
\end{align*}
 We then apply Lemmas \ref{clem1} and \ref{clem2} to the relation
\begin{align*}
\|\hat v_{j,n}\|^2-1&=\| \tilde v_{j,n}\|^2-2(x_n^\top \tilde v_{j,n})^2+(x_n^\top \tilde v_{j,n})^2 -1\\
&=\| \tilde v_{j,n}\|^2-1-(x_n^\top \tilde v_{j,n})^2 
\end{align*}
to find
\begin{align*}
 \big|\|\hat v_{j,n}\|^2-1\big|&\leq\big|\| \tilde v_{j,n}\|^2-1\big|+\big|(x_n^\top \tilde v_{j,n})^2\big|\\
 &\leq Q\tau\sum_{l=1}^{j-1}\|\hat v_{l,n}-v_{l,n}\|+Q\tau^2+Q\tau^4\\
 &\leq Q\tau\sum_{l=1}^{j-1}\|\hat v_{l,n}-v_{l,n}\|+Q\tau^2,
\end{align*}
which completes the proof.
\end{proof}

\section{Numerical analysis for semi-implicit scheme}\label{sec4}
We prove error estimate for the semi-implicit scheme (\ref{cFDsadk}) by performing a multi-variable circulating induction procedure to gradually decouple the quantities of interest.

\subsection{Quantification of $\tilde v_{i,n}-v_{i,n}$}
For $\bar G>Q_3Q_4+kQ_1$ where $Q_1$ and $Q_3$ are introduced in Lemma \ref{clem3} and $Q_4>1$ represents the bound of $\{\tilde v_{j,n}\}_{j=1,n=0}^{k,N}$ (cf. (\ref{xbndmh})), there exists an intermediate constant $G>0$ such that
$$\bar G>Q_3Q_4+kG\text{ and } G>Q_1. $$
In particular, as $Q_4>1$, we have $\bar G>Q_3$. Then for $\tau$ small enough the following inequalities hold:
\begin{equation}\label{GG} \begin{array}{c}
\ds\frac{Q_0\tau k\bar G+Q_1+kG^2\tau^2}{(1-Q_2\tau^3 k\bar G-Q_3\tau^2-kG^2\tau^4)^{1/2}}\leq G,\\[0.15in]
\ds \frac{Q_4(Q_2\tau k \bar G+Q_3+kG^2\tau^2)+kG}{(1-Q_2\tau^3 k\bar G-Q_3\tau^2-kG^2\tau^4)^{1/2}}\leq \bar G.
\end{array}  \end{equation}
In subsequent proofs, we always choose sufficiently small step size $\tau$ such that the condition (\ref{GG}) is satisfied.

\begin{thm}\label{lem2k}
Under the condition (\ref{GG}), the following estimate holds for $1\leq n\leq N$:
$$\|v_{i,n}-\hat v_{i,n}\|\leq \bar G\tau^2,~~1\leq i\leq k.$$
\end{thm}
\begin{rem}
The $\tilde v_{i,n}$ on the left-hand side of the third equation of (\ref{cFDsadk}) could be split as
$$\tilde v_{i,n}=v_{i,n}-(v_{i,n}-\hat v_{i,n})-(\hat v_{i,n}-\tilde v_{i,n}), $$
where the last two right-hand side terms are $O(\tau^2)$ terms according to Lemma \ref{clem1} and this theorem. Thus we reach the following relation that plays a key role in error estimates:
\begin{equation}\label{vvsplit}
\tilde v_{i,n}=v_{i,n}+O(\tau^2).
 \end{equation}
\end{rem}
\begin{proof}
We prove this theorem by induction for the following two relations:
\begin{equation*}
\begin{array}{l}
\ds (\mathbb A):~ \max_{m<i\leq k}\|\hat v_{i,n}^\top v_{m,n}\|\leq G\tau^2 \text{ for some }1\leq m\leq k-1;\\[0.05in]
\ds(\mathbb B):~\|v_{j,n}-\hat v_{j,n}\|\leq \bar G\tau^2\text{ for some }1\leq j\leq k.
\end{array}
\end{equation*}
We first declare that if
 \begin{equation}\label{induc}
 (\mathbb A) \text{ holds for }1\leq m\leq m^*-1\text{ and }(\mathbb B)\text{ holds for }1\leq j\leq m^*
 \end{equation}
  for some $1\leq m^*<k-1$, then
\begin{equation}\label{dec1}
 (\mathbb A) \text{ holds for } m= m^*\text{ and }(\mathbb B)\text{ holds for } j= m^*+1.
 \end{equation}
To show this, we apply Lemma \ref{clem3} and the induction hypotheses (\ref{induc}) to bound $Y_{m^*,n}$ by
\begin{equation}\label{bnd0Y}
 \begin{array}{rl}
 \ds Y_{m^*,n}&\ds=\bigg(\|\hat v_{m^*,n}\|^2-\sum_{j=1}^{m^*-1}(\hat v_{m^*,n}^\top v_{j,n})^2\bigg)^{1/2}\\
& \ds\in \bigg[1\pm \bigg(Q_2\tau\sum_{l=1}^{m^*-1}\|\hat v_{l,n}-v_{l,n}\|+Q_3\tau^2+(m^*-1) G^2\tau^4\bigg)\bigg]^{1/2} \\[0.15in]
& \ds\in \big[1\pm \big(Q_2 (m^*-1)\bar G\tau^3+Q_3\tau^2+(m^*-1) G^2\tau^4\big)\big]^{1/2} .
 \end{array}
\end{equation}
 We then invoke the induction hypotheses (\ref{induc}), (\ref{bnd0Y}), the condition (\ref{GG}) and Lemma \ref{clem3} into the expression of $\hat v_{i,n}^\top v_{m^*,n}$ to obtain for $m^*<i\leq k$
$$\begin{array}{rl}
\ds | \hat v_{i,n}^\top v_{m^*,n}|&\ds=\frac{1}{Y_{m^*,n}}\bigg|\hat v_{i,n}^\top\hat v_{m^*,n}-\sum_{j=1}^{m^*-1}(\hat v_{m^*,n}^\top v_{j,n})(\hat v_{i,n}^\top v_{j,n})\bigg| \\[0.15in]
&\ds\leq\frac{1}{Y_{m^*,n}}\bigg( Q_0\tau\sum_{l=1}^{m^*}\|\hat v_{l,n}-v_{l,n}\|+Q_1\tau^2+(m^*-1)G^2\tau^4\bigg) \\[0.2in]
&\ds \leq \frac{Q_0\tau m^*\bar G+Q_1+(m^*-1)G^2\tau^2}{(1-Q_2\tau^3 (m^*-1)\bar G-Q_3\tau^2-(m^*-1)G^2\tau^4)^{1/2}}\tau^2\leq G\tau^2,
\end{array}$$
which implies that $(\mathbb A)$ holds for $m=m^*$. We then use Lemma \ref{clem3} and $(\mathbb A)$ with $1\leq m\leq m^*$ to bound $Y_{m^*+1,n}$ in an analogous manner as (\ref{bnd0Y}):
 \begin{equation}\label{bndY}
 \begin{array}{rl}
 \ds Y_{m^*+1,n}\in \big[1\pm \big(Q_2 m^*\bar G\tau^3+Q_3\tau^2+m^* G^2\tau^4\big)\big]^{1/2} ,
 \end{array}
\end{equation}
which implies
$$|1-Y_{m^*+1,n}|\leq |1-Y_{m^*+1,n}^2|\leq Q_2 m^*\bar G\tau^3+Q_3\tau^2+m^* G^2\tau^4. $$
We invoke this and $(\mathbb A)$ with $1\leq m\leq m^*$ in $v_{m^*+1,n}-\hat v_{m^*+1,n}$ to get
\begin{equation}\label{diff}
\begin{array}{l}
\ds \|v_{m^*+1,n}-\hat v_{m^*+1,n}\|\\[0.05in]
\ds\qquad=\frac{1}{Y_{m^*+1,n}}\bigg\|(1-Y_{m^*+1,n})\hat v_{m^*+1,n}-\sum_{j=1}^{m^*}(\hat v_{m^*+1,n}^\top v_{j,n})v_{j,n}\bigg\|\\[0.2in]
\ds\qquad\leq \frac{Q_4(Q_2\tau m^* \bar G+Q_3+m^*G^2\tau^2)+m^*G}{(1-Q_2\tau^3 m^*\bar G-Q_3\tau^2-m^*G^2\tau^4)^{1/2}}\tau^2\leq \bar G\tau^2,
\end{array}
\end{equation}
which implies that $(\mathbb B)$ holds for $j=m^*+1$. Therefore, the declaration (\ref{induc})-(\ref{dec1}) is correct and we remain to show that $(\mathbb A)$ holds for $m=1$ and $(\mathbb B)$ holds for $1\leq j\leq 2$ in order to start the mathematical induction. We apply Lemma \ref{clem3} to obtain
 $$\|\hat v_{1,n}-v_{1,n}\|=\bigg\|\frac{\hat v_{1,n}}{\|\hat v_{1,n}\|}(\|\hat v_{1,n}\|-1)\bigg\|\leq ||\hat v_{1,n}\|^2-1|\leq Q_3\tau^2\leq \bar G\tau^2, $$
 which is the relation $(\mathbb B)$ with $j=1$. Based on this, $(\mathbb A)$ with $m=1$ and $(\mathbb B)$ with $j= 2$ can be proved following exactly the same procedure as (\ref{bnd0Y})-(\ref{diff}), which completes the proof.
 \end{proof}

\subsection{Error estimate}
We prove error estimates for the semi-implicit scheme (\ref{cFDsadk}) of sphere-constrained high-index saddle dynamics (\ref{csadk}) by analyzing the following errors:
$$e^x_n:=x(t_n)-x_n,~~e^{v_i}_n:=v_i(t_n)-v_{i,n},~~1\leq n\leq N,~~1\leq i\leq k.$$

\begin{thm}\label{thmevk}
Under the Assumption A, the following estimate holds for the semi-implicit scheme (\ref{cFDsadk}) for $\tau$ sufficiently small:
$$\max_{1\leq n\leq N}\{\|e^x_{n}\|, \|e^{v_1}_n\|,\cdots,\|e^{v_k}_n\|\}\leq Q\tau,~~1\leq n\leq N. $$
Here $Q$ is independent from $\tau$, $n$ and $N$.
\end{thm}
\begin{proof}
To bound $e^x_n$, we derive the reference equation from the first equation of (\ref{csadk}) via the forward Euler discretization
\begin{align*}
x(t_{n}) &=x(t_{n-1})+\tau\bigg(I -x(t_{n-1})x(t_{n-1})^\top\\
&\qquad\qquad-2\sum_{j=1}^k v_j(t_{n-1})v_j(t_{n-1})^\top \bigg)F(x(t_{n-1}))+O(\tau^2). 
\end{align*}
We then apply (\ref{xtx}) and (\ref{xxmh}) to reformulate (\ref{lalala}) as
\begin{align*}
 x_{n} &=x_{n-1}+(x_{n}-\tilde x_n)\\
 &~~+\tau\bigg(I- x_{n-1}x_{n-1}^\top-2\sum_{j=1}^k v_{j,n-1}v_{j,n-1}^\top \bigg)(\mathcal L \tilde x_{n}+\mathcal N(x_{n-1}))+O(\tau^2).\nonumber\\
 &=x_{n-1}+(x_{n}-\tilde x_n)\\
 &~~+\tau\bigg(I- x_{n-1}x_{n-1}^\top-2\sum_{j=1}^k v_{j,n-1}v_{j,n-1}^\top \bigg)F(x_{n-1})\\
 &~~+\tau\bigg(I- x_{n-1}x_{n-1}^\top-2\sum_{j=1}^k v_{j,n-1}v_{j,n-1}^\top \bigg)\mathcal L(\tilde x_n-x_{n-1})+O(\tau^2).\nonumber\\
 &=x_{n-1}+\tau\bigg(I- x_{n-1}x_{n-1}^\top-2\sum_{j=1}^k v_{j,n-1}v_{j,n-1}^\top \bigg)F(x_{n-1})+O(\tau^2).
\end{align*}
By this means, the original semi-implicit scheme of $x$ in (\ref{cFDsadk}) is converted to the explicit scheme to facilitate the analysis. Based on the above two equations, we follow the same derivations in \cite[Theorem 4.2]{Z3c} to obtain
\begin{equation}\label{xref}
\|e^x_n\|\leq Q\tau\sum_{m=1}^{n-1}\sum_{j=1}^k \|e^{v_j}_m\|+Q\tau,~~1\leq n\leq N.
\end{equation}

To estimate $e^{v_i}_n$, we derive the reference equation from the third equation of (\ref{csadk}) via the backward Euler discretization for $1\leq i\leq k$:
\begin{equation*}
\begin{array}{l}
\ds v_i(t_{n})=v_i(t_{n-1})+\tau\bigg( I- x(t_n)x(t_n)^\top-v_i(t_{n})v_i(t_{n})^\top\\
\ds\hspace{0.5in}-2\sum_{j=1}^{i-1}v_j(t_{n})v_j(t_{n})^\top\bigg)J(x(t_{n}))v_i(t_{n})+ \tau x(t_n)v_i(t_n)^\top F(x(t_n))+O(\tau^2)\\
\ds\qquad~\,=v_i(t_{n-1})+\tau\bigg( I- x(t_n)x(t_n)^\top-2\sum_{j=1}^{i-1}v_j(t_{n})v_j(t_{n})^\top\bigg)J(x(t_{n}))v_i(t_{n})\\
\ds\hspace{0.5in}-\tau v_i(t_{n-1})v_i(t_{n-1})^\top J(x(t_n))v_i(t_{n-1})\\[0.05in]
\ds\hspace{0.5in}+\tau x(t_n)v_i(t_n)^\top F(x(t_n))+O(\tau^2)+\mathcal A_n
\end{array}
\end{equation*}
where
$$\begin{array}{l}
\ds \mathcal A_n=\tau\big(v_i(t_{n})v_i(t_{n})^\top J(x(t_n))v_i(t_{n})\\[0.05in]
\ds\qquad\qquad-v_i(t_{n-1})v_i(t_{n-1})^\top J(x(t_n))v_i(t_{n-1})\big)=O(\tau^2).
\end{array}
$$
We then apply (\ref{vvsplit}) to rewrite the third scheme of (\ref{cFDsadk}) as
\begin{align*}
 v_{i,n}&=v_{i,n-1}+(v_{i,n}-\tilde v_{i,n})\\
 &\qquad+\tau\bigg( I-x_{n}x_{n}^\top-2\sum_{j=1}^{i-1}v_{j,n}v_{j,n}^\top\bigg)J(x_{n})( v_{i,n}+O(\tau^2))\\
&\qquad-\tau v_{i,n-1}v_{i,n-1}^\top J(x_{n})v_{i,n-1}+\tau x_{n}(v_{i,n}+O(\tau^2))^\top F(x_{n})\\
&=v_{i,n-1}+\tau\bigg( I-x_{n}x_{n}^\top-2\sum_{j=1}^{i-1}v_{j,n}v_{j,n}^\top\bigg)J(x_{n})v_{i,n}\\
&\qquad-\tau v_{i,n-1}v_{i,n-1}^\top J(x_{n})v_{i,n-1}+\tau x_{n}v_{i,n}^\top F(x_{n})+O(\tau^2).
\end{align*}
Based on the above two equations, we follow almost the same derivations as \cite[Theorem 4.2]{Z3c} to derive the estimate of $e^{v_i}_n$ as
$$ \sum_{i=1}^k\|e^{v_i}_n\|\leq Q\tau,$$
and we invoke this in (\ref{xref}) to complete the proof.
\end{proof}

\section{Numerical experiments}
We carry out a simple numerical experiment to test the convergence rate (denoted by CR) of the scheme (\ref{cFDsadk}). A detailed comparison between semi-implicit and explicit methods for unconstrained high-index saddle dynamics could be found in \cite{LuoSI}, which has already indicated the advantages of the semi-implicit method.
We apply the Rosenbrock type function
\begin{equation*} E(x_1,x_2,x_3)=a(\sqrt{3}x_2-3x_1^2)^2+b(\sqrt{3}x_1-1)^2+a(\sqrt{3}x_3-3x_2^2)^2+b(\sqrt{3}x_2-1)^2.
\end{equation*}
 For $(a,b)=(-1,5.5)$, the point 
$$x_*=\mathcal N(1,1,1)=\frac{1}{\sqrt{3}}(1,1,1)$$
is an index-1 saddle point of the Rosenbrock type function, while for  $(a,b)=(-0.5,1.5)$, $x_*$  is an index-2 saddle point. 
 We apply the semi-implicit scheme (\ref{cFDsadk}) to compute the saddle points for these two cases under
 $T=10$ and different initial conditions
\begin{align*}
&\text{(a) } x_0=\mathcal N(0.8,1,1),~~v_{1,0}=\mathcal N(1,-0.4,-0.4);\\
&\text{(b) } x_0=\mathcal N(1,1,1.4),~~v_{1,0}=\mathcal N(-1,1,0);\\
&\text{(c) } x_0=\mathcal N(0.8,1,1),~~v_{1,0}=\mathcal N(1,-0.4,-0.4),~~v_{2,0}=\mathcal N(0,1,-1);\\
&\text{(d) } x_0=\mathcal N(1,1,1.4),~~v_{1,0}=\mathcal N(-1,1,0),~~v_{2,0}=\mathcal N(-0.7,-0.7,1).
\end{align*}
As the exact trajectory of the constrained high-index saddle dynamics (\ref{csadk}) is in general not available, we use the numerical solution computed under $\tau=2^{-13}$ to serve as the reference solution.  Numerical results are presented in Tables \ref{table1}--\ref{table4}, which indicates the first-order accuracy of the semi-implicit scheme (\ref{cFDsadk}) as proved in Theorem \ref{thmevk}.
\begin{table}[h!]
\setlength{\abovecaptionskip}{0pt}
\centering
\caption{CR of computing the index-1 saddle point under the initial condition (a).}
\begin{tabular}{ccccc} \cline{1-5}
$\tau$& $\max_n \|e^x_n\|$ & CR &  $\max_n \|e^{v_1}_n\|$ &CR\\
\cline{1-5}		
$2^{-6}$&	1.65E-02	&	&9.95E-02	&\\
$2^{-7}$&	8.29E-03&	0.99 &	4.47E-02	&1.16\\ 
$2^{-8}$&	4.09E-03& 1.02 	&2.08E-02	&1.10 \\
$2^{-9}$&	1.98E-03&	1.04 	&9.83E-03&	1.08\\ 
				\hline
			\end{tabular}
			\label{table1}
		\end{table}

\begin{table}[h!]
\setlength{\abovecaptionskip}{0pt}
\centering
\caption{CR of computing the index-1 saddle point under the initial condition (b).}
\begin{tabular}{ccccc} \cline{1-5}
$\tau$& $\max_n \|e^x_n\|$ & CR &  $\max_n \|e^{v_1}_n\|$ &CR\\
\cline{1-5}		
$2^{-6}$&	1.03E-02&		&2.02E-02	&\\
$2^{-7}$&	4.84E-03&	1.09 &	9.53E-03	&1.08\\ 
$2^{-8}$&	2.32E-03&	1.06 &	4.59E-03	&1.05 \\
$2^{-9}$&	1.11E-03&	1.06 &	2.20E-03&	1.06 \\
				\hline
			\end{tabular}
			\label{table2}
		\end{table}

\begin{table}[h!]
\setlength{\abovecaptionskip}{0pt}
\centering
\caption{CR of computing the index-2 saddle point under the initial condition (c).}
\begin{tabular}{ccccccc} \cline{1-7}
$\tau$& $\max_n \|e^x_n\|$ & CR &  $\max_n \|e^{v_1}_n\|$ &CR&  $\max_n \|e^{v_2}_n\|$ &CR\\
\cline{1-7}		
$2^{-6}$&	1.67E-03&		&6.06E-02&		&6.06E-02&	\\
$2^{-7}$&	7.90E-04&	1.08& 	2.87E-02&	1.08 &	2.87E-02&	1.08 \\
$2^{-8}$&	3.80E-04&	1.05& 	1.38E-02&	1.06 &	1.38E-02&	1.06 \\
$2^{-9}$&	1.82E-04&	1.06& 	6.60E-03&	1.06& 	6.60E-03&	1.06 \\
				\hline
			\end{tabular}
			\label{table3}
		\end{table}

\begin{table}[h!]
\setlength{\abovecaptionskip}{0pt}
\centering
\caption{CR of computing the index-2 saddle point under the initial condition (d).}
\begin{tabular}{ccccccc} \cline{1-7}
$\tau$& $\max_n \|e^x_n\|$ & CR &  $\max_n \|e^{v_1}_n\|$ &CR&  $\max_n \|e^{v_2}_n\|$ &CR\\
\cline{1-7}		
$2^{-6}$&	2.65E-03&		&3.53E-02	&	&3.52E-02	&\\
$2^{-7}$&	1.28E-03&	1.05 &	1.69E-02&	1.06& 	1.69E-02&	1.06 \\
$2^{-8}$&	6.22E-04&	1.04 &	8.21E-03&	1.05 &	8.18E-03&	1.05 \\
$2^{-9}$&	2.99E-04&	1.06 &	3.94E-03&	1.06 &	3.93E-03&	1.06 \\
				\hline
			\end{tabular}
			\label{table4}
		\end{table}

\section{Concluding remarks}
In this paper we prove error estimates for the semi-implicit numerical scheme of sphere-constrained high-index saddle dynamics, which ensures the accuracy of performing the saddle dynamics in finding saddle points and constructing the solution landscape for constrained problems. The main difficulties we overcome lie in the semi-implicit treatment on the schemes and the coupling among the dynamics, the retraction, the vector transport and the orthonormalization procedure. Numerical experiments are performed to substantiate the theoretical findings.

There are potential extensions of the current work that deserve further exploration. For instance, the dimer method \cite{Dimer} could be used in (\ref{csadk}) to approximate the product of the Hessian matrix and the vector for efficient computation and storage, which leads to the shrinking-dimer sphere-constrained high-index saddle dynamics as the unconstrained case \cite{ZhaCSIAM}. Then the semi-implicit method could be applied to improve the numerical stability that remains to be analyzed.

Furthermore, the ideas and techniques could be employed and improved to analyze the semi-implicit numerical scheme for high-index saddle dynamics constrained by $m$ equalities \cite[Equation 24]{CHiSD2021}:
\begin{equation}\label{ccsd}
\left\{
\begin{array}{l}
\ds \frac{dx}{dt} =\bigg(I -2\sum_{j=1}^k v_jv_j^\top \bigg)F(x),\\[0.075in]
\ds \frac{dv_i}{dt}=\bigg(I-v_iv_i^\top-2\sum_{j=1}^{i-1}v_jv_j^\top\bigg)\mathcal H(x)[v_i] \\
\ds\qquad\qquad-A(x)\big(A(x)^\top A(x)\big)^{-1}\bigg(\nabla^2c(x)\frac{dx}{dt}\bigg)^\top v_i,~~1\leq i\leq k.
\end{array}
\right.
\end{equation}
Here $c(x)=(c_1(x),\cdots,c_m(x))=0$ represents the $m$ equality constraints and 
$$A(x)=(\nabla c_1(x),\cdots,\nabla c_m(x)).$$
 The sphere-constrained high-index saddle dynamics (\ref{csadk}) is a special case of (\ref{ccsd}) with one equality constraint
$$c_1(x)=\|x\|-1=0. $$
In the generalized constrained saddle dynamics (\ref{ccsd}), $\mathcal H(x)$ refers to the Riemannian Hessian \cite{CHiSD2021}, which is difficult to compute and approximate in practice that we will investigate in the near future.

\bigskip


\end{document}